\documentclass[12pt]{article}
\usepackage[russian,english]{babel}
\usepackage[cp1251]{inputenc}
\usepackage{srcltx}
\usepackage{amssymb,amsmath,amsfonts,amsthm}
\usepackage{cite}
\newtheorem{lem}{Lemma}

\newtheorem{theorem}{Theorem}

\begin{document}

UDC 512.542

\begin{center}\textbf{On Thompson's conjecture for alternating and symmetric groups }\end{center}

\begin{center}I.B. Gorshkov \footnote{The work is supported by  Russian Science Foundation (project
14-21-00065)}\end{center}

\textbf{1. Introduction}
\medskip

In this paper, all groups are finite.
Let $G$ be a finite group. If $g\in G$ that $g^G$ denotes the conjugacy class in $G$ containing $g$ and $|g^G|$ its size, $C_G(g)$ denotes the centralizer of $g$ in $G$. We use $N(G)$ for the set $\{n| \exists g\in G$ such that $|g^G|=n\}$. Thompson in
1987 posted the following conjecture with respect to $N(G)$.
\medskip

\textbf{Thompson's Conjecture (see \cite{Kour}, Question 12.38)}. {\it If $L$ is a finite
simple non Abelian group, $G$ is a finite group with trivial
center, and $N(G)=N(L)$, then $G\simeq L$.}

\medskip

Thompson’s conjecture was proved valid for many finite simple groups of Lie type.
Let $Alt_n$ be the alternating group of degree $n$, $Sym_n$ be the symmetric group of degree $n$. Alavi and Daneshkhah
proved that the groups $Alt_n$ with $n=p$, $n=p+1$, $n=p+2$ and
$p$ prime are characterized by $N(G)$ (see \cite{AD}). Recently, the groups $Alt_{10}, Alt_{16}$,
and $Alt_{22}$ were proved valid for this conjecture (see \cite{VasT}, \cite{Gor}, \cite{Xu}).

Our main result is the following.

\begin{theorem}
If $G$ is a finite group such that $N(G)=N(Alt_n), n\geq5$ or $N(G)=N(Sym_n)$, $n\geq23$, then $G$ is non-solvable.
\end{theorem}

\textbf{ 2. Notation and preliminary results}

\medskip

Given a finite group $G$, denote by $\pi(G)$ the set of prime divisors of its order.
\begin{lem}\label{pi}
Let $G$ and $H$ are finite groups, center of $G$ is trivial and $N(G) = N(H)$.
Then $\pi(G)\subseteq \pi(H)$.
\end{lem}
\textsl{Proof.}  In the same way as \cite[Lemma 3]{VasT}.

\begin{lem}[{\rm \cite[Lemma 5]{VasT}}]\label{factor}
Let $K$ be a normal subgroup of a finite group $G$, and $\overline{G}=G/K$. If $\overline{x}$ is the image of an element $x$ of $G$ to $\overline{G}$, them $|\overline{x}^{\overline{G}}|$ divides $|x^G|$. Moreover, if $(|x|,|K|)=1$, then $C_{\overline{G}}(\overline{x})=C_G(x)K/K$.
\end{lem}

\begin{lem}\label{centr}
Let $x,y\in G, (|x|,|y|)=1$, $xy=yx$. Then $C_G(xy)=C_G(x)\cap
C_G(y)$.
\end{lem}
\textsl{Proof.}
The proof is trivial.

\begin{lem}[{\rm \cite[Theorem 4.34]{Isaacs}}]\label{Fitting}
Let $A$ act via automorphisms on an Abelian group $G$, and suppose that $(|G|,|A|)=1$. Then $G=C_G(A)\times[G,A]$.
\end{lem}

\begin{lem}\label{Frobenius}
Let $g$ act via automorphisms on an Abelian group $G$, and suppose that $(|G|,|g|)=1$. Then $|g|$ divides $|[G,g]|-1$.
\end{lem}

\begin{lem}\label{vas}
If $G$ is finite group and there exists a prime $p\in \pi(G)$ such that $p^2$ does not divide $|x^G|$ for all $x$ in $G$. Then a Sylow $p$-subgroup of $G$ is Abelian.
\end{lem}
\textsl{Proof.} In the same way as \cite[Lemma 4]{VasT}.

Let $\Theta\subseteq \mathbb{N}, |\Theta|<\infty$, $\Gamma(\Theta)$ is directed graph, the set vertices is equal to $\Theta$  and $\overrightarrow{ab}$ if $a$ divides $b$, $h(\Theta)$ is length of a maximal path in the graph $\Gamma(\Theta)$. Let $\Gamma(G)=\Gamma(N(G))$, $V_i\in\{Alt_i, Sym_i\}$, $G$ be a finite solvable group such that $N(G)=N(V_n)$ where $n\geq5$ if $V_n\simeq Alt_n$ and $n\geq 23$ if $V_n\simeq Sym_n$, $\Omega=\{t| n/2<t\leq n, t -$ prime $\}$, $p$ be maximal number of $\Omega$. From Lemma \ref{pi} we obtain that $\pi(G)\supseteq\pi(V_n)$ in particular $\Omega \subseteq \pi(G)$.

\begin{lem}\label{omega}
If $n>1361$, then $|\Omega|>\log_2(n!/(p)!)$,
\end{lem}
\textsl{Proof.} Let $\phi (x)$ be the number of primes less or equal $x$,
From \cite[Section 35, \S 1]{Buh} it follows that $0,921\cdot (x/\ln(x))<
\phi(x) < 1,106\cdot (x/\ln(x))$ for $10<x$.
Hence $|\Omega|=\phi(n)-\phi(n/2)\geq \frac{0,921\cdot n}{\ln(n)}- \frac{1,106\cdot
n/2}{\ln(n/2)}$. It is proved in \cite{Baker} that $n-p<n^{0,525}$.
Using this assertions we get that lemma is valid if $n\geq1000000$. If $1361\leq n<1000000$ then using \cite{GAP} we obtain our assertion. The lemma is proved.

\begin{lem}\label{ConClass}
Let $t\in \Omega$, $\alpha\in N(G)$ and $t\not \in \pi(\alpha)$. Then $\alpha=|V_n|/(t|C|)$ or $\alpha=|V_n|/(|V_{t+i}||B|)$ where $C=C_{V_{n-t}}(g)$, for any $g\in V_{n-t}$, $t+i\leq n$, $B=C_{V_{n-t-i}}(h), h\in V_{n-t-i}$ and $h$ moves $n-t-i$ points.
\end{lem}
\textsl{Proof.}
The proof is trivial.
\medskip

 Let $\Phi_t=\{\alpha\in N(G)| \alpha=|V_n|/(t|C|)$, $C=C_{V_{n-t}}(g)$, for any $g\in V_{n-t}\}$, $\Psi_t=\{\alpha\in N(L)| \alpha=|V_n|/(|V_{t+i}||B|$, $i\geq0, t+i<n-1$, $B=C_{V_{n-t-i}(g)}, g\in V_{n-t-i}$ and $g$ moves $n-t-i$ points $\}$. Let us remark that the definitions of sets $\Phi$ and $\Psi$ do not imply that $t$ is prime.

\begin{lem}\label{hz}
\begin{enumerate}
\item{If $n-t=2$, then $h(\Psi_t)\leq 1$ }
\item{If $n-t=3$, then $h(\Psi_t)\leq 2$ }
\item{If $n-t=4$, then $h(\Psi_t)\leq 3$ }
\item{If $n-t=5$, then $h(\Psi_t)\leq 5$ }
\item{If $n-t=6$, then $h(\Psi_t)\leq 6$ }
\item{If $n-t=7$, then $h(\Psi_t)\leq 8$ }
\item{If $n-t=8$, then $h(\Psi_t)\leq 11$ }
\item{If $n-t=9$, then $h(\Psi_t)\leq 14$ }
\item{If $n-t=10$, then $h(\Psi_t)\leq 18$ }
\item{If $n-t=11$, then $h(\Psi_t)\leq 21$ }
\item{If $n-t=12$, then $h(\Psi_t)\leq 26$ }
\item{If $n-t=13$, then $h(\Psi_t)\leq 30$ }
\item{If $n-t=18$, then $h(\Psi_t)\leq 69$ }

\end{enumerate}
\end{lem}\textsl{Proof.} By definition we have $h(\Psi_t)\leq \sum_{i\in\{1,..,n-t\}}h(R_i)$ where $R_i=\{|g^{V_i}|, g\in V_i,$ and $g$ moved $i$ point $\}$.
Let us use \cite{GAP} for calculation $R_i$. The lemma is proved.

\begin{lem}\label{abelian}
If $G$ is solvable group, then Hall $\Omega$-subgroup is an Abelian.
\end{lem}\textsl{Proof.}
It follows from Lemma \ref{vas} that Sylow $t$-subgroup of $G$ is Abelian for any  $t\in \Omega$. Assume that Hall $\Omega$-subgroup of $G$ is not Abelian. It follows that there exists $\{t_i, t_j\}\subset\Omega$ such that a Hall $\{t_i,t_j\}$-subgroup $H$ of $G$ is non Abelian. Let $R<G$ be a maximal normal subgroup such that $\overline{H}=HR/R$ is non Abelian. Hence in group $\overline{G}=G/R$ there exist normal $t$-subgroup $T$, $t\in\{t_1, t_2\}$, such that there exist $g\in \overline{G}, |g|\in\{t_1,t_2\}\setminus \{t\}$ and $[T,g]\neq 1$. By Lemma \ref{Fitting}, $T=C_T(g)\times[T,g]$. Using Lemma \ref{Frobenius} we get that $|g|$ divides $[T,g]-1$. But $|g|$ not divides $t-1$, therefore $|[T,g]|>t$. From Lemma \ref{factor} it follows that $|g^{\overline{G}}|_t>t$ and there exists $n\in N(G)$ such that $t^2$ divides $n$; a contradiction. The lemma is proved.

\begin{lem}\label{abelian2}
Let $G$ by solvable. If there exists $r$ prime such that $p+1<2r\leq n$ then a Hall $\{r\}\cup\Omega$-subgroup is isomorphic direct product of $T\times R$ where $T$ is a Hall $\Omega$-subgroup and $R$ is a Sylow $r$-subgroup. In particular, if $g\in T$ then $|g^G|$ is not divisible by $r$.
\end{lem}\textsl{Proof.} The proof is analogous to that of Lemma \ref{abelian}.

\begin{lem}\label{solvable}
Let $T$ be a Hall $\Omega$-subgroup of $G$, $\Theta=\{|g^G|$ for all $g\in T\}$.
If $|\Omega|>h(\Theta)$, then $G$ are non-solvable.
\end{lem}
\textsl{Proof.} Assume that $G$ is solvable and $|\Omega|>h(\Theta)$. From Lemma \ref{abelian} it follows that a Hall $\Omega$-subgroup of $G$ is Abelian. Let $g_1\in T, |g_1|=t_1, t_1\in \Omega$ is minimal of $\Omega$. It follows that there exist $r_1\in C_G(g_1)$ such that intersection $C_G(r_1)\cap T$ is Hall $\Omega$-subgroup of $C_G(r_1)$, and $t_2$ divides $|r_1^G|$ where $t_2$ is minimal of $\Omega\setminus\{t_1\}$. Hence there exist $g_2\in T, |g_2|=t_2$ such that $C_G(g_1)\neq C_G(g_2)$. Since $T$ is Abelian from Lemma \ref{centr} it follows that $|(g_1g_2)^G|>|g_1^G|$ and $|g_1^G|$ divides $|(g_1g_2)^G|$. Repeating this procedure $|\Omega|$ times we obtain a set $\Sigma=\{g_1, g_1g_2, g_1g_2g_3,..., g_1g_2...g_{|\Omega|}\}$ and $|g_1^G| | |(g_1g_2)^G| |...||g_1g_2...g_{|\Omega|}|$. Thus $h(\Theta)\geq |\Omega|$; a contradiction. The lemma is proved.

\medskip

\textbf{3. Proof of the main theorem}

\medskip

\begin{lem}
If $n>1361$ then $G$ is non-solvable.
\end{lem}
\textsl{Proof.} Assume that $G$ is solvable group. Let $T$ be a Hall $\Omega$-subgroup of $G$. Using Lemma \ref{ConClass} we get that $|g^G|\in\Psi_p$ for any $g\in T$. Let us show that $|\Omega|>h(\Psi_p)$. Let $h_1,...,h_k\in \Psi_p$ such that $h_1|h_2|...|h_k$. Therefore $2h_1\leq h_2, 2h_2\leq h_3...2h_{k-1}\leq h_k$. If $l\in\Psi_p$ then $l|n!/(p)!$. Hence $h(\Psi_p)\leq log_2(n!/(p)!)$. Using Lemma \ref{omega} we get $|\Omega|>log_2(n!/(p)!\geq h(\Psi_p)$. From Lemma \ref{solvable} if follow that $G$ is non-solvable; a contradiction. The lemma is proved.

\begin{lem}
If $22<n\leq1361$ then $G$ is non-solvable.
\end{lem}
\textsl{Proof.} Assume that $G$ is solvable group. Let $T$ be a Hall $\Omega$-subgroup of $G$. Using Lemma \ref{ConClass} we get that $|g^G|\in\Psi_p$ for any $g\in T$. From Lemma \ref{abelian2} it follows that $|g^G|$ is not divisible by $r$ where $g\in T$ and $p+1<2r\leq n$, $r$ is prime. Hence $|g^G|=|V_n|/(V_{2r+i}|C|)\in\Psi_{2r}$ where $0\leq i \leq n-2r+2$, and $C=C_{V_{n-2r-i}(g)}, g\in V_{n-t-i}, 0\leq i\leq n-2r+2$. Let $\Theta=\{|g^G|$ for all $g\in T\}$. Combining Lemma \ref{hz}, \cite{GAP}, and Lemma \ref{solvable}, we get a contradiction. The lemma is proved.

\begin{lem}
If $5\leq n\leq22$ then $G$ is non-solvable.
\end{lem}
\textsl{Proof.} The assertion follows from \cite{AD}, \cite{VasT} and \cite{Gor}.

\thispagestyle{empty}
I.B. Gorshkov\\
Sobolev Institute of Mathematics,\\
4 Acad. Koptyug avenue, 630090, Novosibirsk, Russia\\
Universidade de Sao Paulo, Brazil\\
E-mail address: ilygor8@gmail.com

\begin{thebibliography}{1}
\bibitem{Kour} V. D. Mazurov, E. I. Khukhro, Eds., The Kourovka Notebook:
Unsolved Problems in Group Theory, Russian Academy of Sciences Siberian Division,
Institute of Mathematics, Novosibirsk, Russia, 18th edition, 2014.

\bibitem{AD} S. H. Alavi and A.Daneshkhah, A new characterization of alternating
and symmetric groups, Journal of Applied Mathematics and Computing, 17(1), pp. 245–258, 2005.

\bibitem{VasT} A. V. Vasil'ev, On Thompson's conjecture, Sibirskie Elektronnye Matematicheskie Izvestiya, vol. 6, pp. 457–464, 2009.

\bibitem{Gor} I. B. Gorshkov, Thompson’s conjecture for simple groups with
a connected prime graph, Algebra and Logic, 51(2), pp. 111–127, 2012.
\bibitem{Xu} M. Xu, Thompson's conjecture for alternating group of degree
22, Frontiers of Mathematics in China, 8(5), pp. 1227–
1236, 2013.

\bibitem{Isaacs} I. M. Isaacs, Finite group theory, Graduate Studies in Mathematics, 92. American Mathematical Society, Providence, RI, 2008. xii+350 pp.

\bibitem{Buh} A. A. Buhshtab,  Number theory // Moskow. 1966

\bibitem{Baker} R. C. Baker, The difference between consecutive primes, II. Proceeding of the London Mathematical Society 83(3), pp. 532-562, 2001.

\bibitem{GAP} The GAP Group, GAP~--- Groups, Algorithms, and Programming, Version 4.4, 2004;
(http://www.gap-system.org).

\end{thebibliography}
\end{document}